\documentclass[a4paper]{amsart}
\date{February 21, 2003}
\usepackage{amsthm}
\usepackage[dvips]{graphicx}
\theoremstyle{plain}
 \newtheorem{theorem}{Theorem}[section]
 
 \newtheorem{lemma}[theorem]{Lemma}
 \newtheorem{corollary}[theorem]{Corollary}
\theoremstyle{remark}
 \newtheorem{definition}[theorem]{Definition}
 
 \newtheorem{example}[theorem]{Example}
 \newtheorem*{remark*}{Remark}
\numberwithin{equation}{section}
\newcommand{\CP}{\mathop{\mathbf CP}}

\newcommand{\bmath}[1]{\mbox{\boldmath $#1$}}
\newcommand{\id}{\mathop{\rm id}}
\newcommand{\Z}{{\bmath{Z}}}
\newcommand{\R}{{\bmath{R}}}
\newcommand{\C}{{\bmath{C}}}
\newcommand{\SL}{\mathop{\rm SL}}

\newcommand{\PSU}{\mathop{\rm PSU}}
\newcommand{\PSL}{\mathop{\rm PSL}}
\newcommand{\Sl}{\mathop{\mathfrak{sl}}}
\newcommand{\N}{\mathop{\mathcal{N}}}
\renewcommand{\L}{\mathop{\mathcal{L}}}
\newcommand{\rank}{\mathop{\rm rank}}
\title[Legendrian curves in $\PSL(2,\C)$]{
   An elementary proof of Small's formula\\
   for null curves in \boldmath{$\PSL(2,\C)$} and\\
   an analogue for Legendrian curves in \boldmath{$\PSL(2,\C)$}
}
\author{Masatoshi Kokubu}
\address[Kokubu]{%
   Department of Natural Science,
   School of Engineering,
   Tokyo Denki University,
   2-2 Kanda-Nishiki-Cho,
   Chiyoda-Ku, Tokyo, 101-8457,
   Japan
}
\email{kokubu@cck.dendai.ac.jp}
\author{Masaaki Umehara}
\address[Umehara]{%
   Department of Mathematics, Graduate School of Science,
   Hiroshima University,
   Higashi-Hiroshima 739-8526, Japan%
}
\email{umehara@math.sci.hiroshima-u.ac.jp}
\author{Kotaro Yamada}
\address[Yamada]{%
   Faculty of Mathematics,
   Kyushu University 36, 
   Higashi-ku, Fukuoka 812-8185, Japan%
}
\email{kotaro@math.kyushu-u.ac.jp}
\dedicatory{
  Dedicated to Professor Katsuei Kenmotsu on his sixtieth birthday
}
\begin{document}
\maketitle
\section{Introduction}
Let $M^2$ be a Riemann surface, 
{\it which might not be simply connected}.
A meromorphic map $F$ from $M^2$ into $\PSL(2,\C)=\SL(2,\C)/\{\pm\id\}$
is a map which is represented as
\begin{equation}\label{eq:meromap}
   F = \begin{pmatrix}
           A & B \\
           C & D
       \end{pmatrix}
        =
       \sqrt{h}\begin{pmatrix}
	   \hat A & \hat B \\
	   \hat C & \hat D
       \end{pmatrix}
       \qquad
        \bigl(AD-BC=1\bigr),
\end{equation}
where $\hat A$, $\hat B$, $\hat C$, $\hat D$
and $h$ are meromorphic functions on $M^2$.
Though $\sqrt{h}$ is a multi-valued function on $M^2$, 
$F$ is well-defined as a $\PSL(2,\C)$-valued mapping.

A meromorphic map $F$ as in \eqref{eq:meromap} is called a 
{\it null curve\/} if the  pull-back of the Killing form by $F$
vanishes, which is equivalent to the condition that the derivative
$F_z=\partial F/\partial z$ with respect to each complex coordinate $z$
is a degenerate matrix everywhere.
It is well-known that the projection of a null curve in $\PSL(2,\C)$
into the hyperbolic $3$-space $H^3=\PSL(2,\C)/\PSU(2)$ gives a constant
mean curvature one surface (see \cite{Br,UY1}).
For a non-constant null curve $F$, we define two meromorphic functions
\begin{equation}\label{eq:G-g}
  G:= \frac{dA}{dC}=\frac{dB}{dD},\qquad
  g:= -\frac{dB}{dA}=-\frac{dD}{dC}.
\end{equation}
(For a precise definition, see Definition~\ref{def:null-gauss} in
Section~\ref{sec:small}).
We call $G$ the {\em hyperbolic Gauss map\/} of $F$ and 
$g$ the {\em secondary Gauss map}, respectively \cite{UY3}.
In 1993, Small \cite{S} discovered the following expression
\begin{equation}\label{eq:small-formula}
 F=\begin{pmatrix}
     G\displaystyle\frac{da}{dG}-a & 
     G\displaystyle\frac{db}{dG}-b \\[6pt]
     \hphantom{G}\displaystyle\frac{da}{dG}\hphantom{-a} & 
     \hphantom{G}\displaystyle\frac{db}{dG}\hphantom{-b}
   \end{pmatrix},\quad
   \left(a:=\sqrt{\frac{dG}{dg}},~
         b:=-g a \right)
\end{equation}
for null curves such that both $G$ and $g$ are non-constant.
({\it We shall give a simple proof of this formula in
      Section~\ref{sec:small}.
      Sa Earp and Toubiana \cite{ET} gave an alternative proof,
      which is quite different from ours. 
      On the other hand, Lima and Roitman \cite{LR}
      explained this formula via the method of Bianchi \cite{Bi} 
      from the 1920's.
      Recently, Small \cite{S1} gave some remarks on this formula
      from the viewpoint of null curves in $\C^4$.})
In this expression, $F$ is expressed by only the derivation of two Gauss
maps.
Accordingly, the formula is valid even if $M^2$ is not simply connected.

By the formula \eqref{eq:small-formula}, it is shown that the set of
non-constant null curves on $M^2$ with non-constant Gauss maps
corresponds bijectively to the set of pairs $(G,g)$ of meromorphic
functions on $M^2$ such that $g\not\equiv a\star G$ 
(that is, $g$ is not identically equal to $a\star G$) for any
$a\in\SL(2,\C)$.
Here, for a matrix $a=(a_{ij})\in\SL(2,\C)$, we denote by $a\star G$ 
the M\"obius transformation of $G$:
\begin{equation}\label{eq:moebius}
    a\star G :=\frac{a_{11}G + a_{12}}{a_{21}G+a_{22}}.
\end{equation}
For this correspondence, see also \cite{UY2}.

On the other hand, according to G\'alvez, Mart\'\i{}nez and Mil\'an
(\cite{GMM1,GMM2}), a meromorphic map
\begin{equation}\label{eq:legendrian-curve}
   E=\begin{pmatrix}
         A & B \\ C & D
     \end{pmatrix}
\end{equation}
from $M^2$ into $\PSL(2,\C)$ is called a {\it Legendrian curve\/} (or a
{\it contact curve}) if the pull-back of the holomorphic contact
form 
\begin{equation}\label{eq:contact}
 D dA -B dC
\end{equation}
on $\PSL(2,\C)$ by $E$ vanishes.
For a Legendrian curve $E$, two meromorphic functions 
\begin{equation}\label{eq:legendrianG}
     G=\frac{A}{C}, \qquad
     G_*=\frac{B}{D}
\end{equation}
are defined. 
In \cite{GMM1}, $G$ and $G_*$ are called the {\it hyperbolic Gauss
maps}.
We define a meromorphic $1$-form $\omega$ on $M^2$ as
\begin{equation}\label{eq:omega}
    \omega := \frac{dA}{B}=\frac{dC}{D}.
\end{equation}
(For a precise definition, see
Definition~\ref{def:invariants-legendrian}
and Lemma~\ref{lem:can-form} in Section~\ref{sec:legendre}.)
We shall call $\omega$ the {\it canonical form}.

As an analogue of the Bryant representation formula \cite{Br,UY1} for
constant mean curvature one surfaces in $H^3$,
G\'alvez, Mart\'\i{}nez and Mil\'an \cite{GMM1} showed that
any simply connected flat surface in hyperbolic $3$-space can be lifted
to a Legendrian curve in $\PSL(2,\C)$, 
where the complex structure of the surface is given so that the second
fundamental form is hermitian.
It is natural to expect that there is a Small-type formula for Legendrian
curves in $\PSL(2,\C)$.

In this paper, we shall give a representation formula for Legendrian 
curves in terms of $G$ and $G_*$ (Theorem~\ref{thm:repG}).
Namely, for an arbitrary pair of non-constant meromorphic functions
$(G,G_*)$ such that $G\not\equiv G_*$ ($G$ is not identically equal to
$G_*$),
the Legendrian curve $E$ with hyperbolic Gauss maps $G$ and $G_*$ is
written as
\begin{equation}\label{eq:repG}
 E = \begin{pmatrix}
      G/\xi&       \xi G_{*}/(G-G_{*}) \\
      1/\xi &      \hphantom{G_*}\xi/(G-G_{*}) 
     \end{pmatrix}\qquad
     \left(\xi = c \exp\int_{z_0}^{z}\frac{dG}{G-G_*}\right),
\end{equation}
where $z_0\in M^2$ be a base point and $c\in\C\setminus\{0\}$ is 
a constant.
As a corollary of this formula, we shall give a Small-type
representation formula for Legendrian curves
(Corollary~\ref{cor:legendre}):
\begin{equation}\label{eq:sm-analogue}
 E=\begin{pmatrix}
    A & dA/\omega \\
    C & dC/\omega 
   \end{pmatrix},\qquad
   \left(C:=i\sqrt{\frac{\omega}{dG}},\quad 
         A:= GC\right).
\end{equation}
It should be remarked that the formula \eqref{eq:sm-analogue}
has appeared implicitly in \cite[page 423]{GMM1} by a different
method.

In Section~\ref{sec:example}, we shall give new examples of flat
surfaces with complete ends using these representation formulas.
Though these examples might have singularities, 
they can be lifted as a Legendrian immersion into the unit cotangent
bundle of $H^3$, and so we call them {\em flat {\rm(}wave{\rm)}
fronts}.
See \cite{KUY}
for a precise definition and global properties of flat fronts with
complete ends.

Small's formula is an analogue of the classical representation
formula for
null curves in $\C^3$, which is closely related to the Weierstrass
representation formula for minimal surfaces in $\R^3$.
For the reader's convenience, we give a simple proof of the
classical formula in the appendix.
\section{A simple proof of Small's formula}
\label{sec:small}
In this section, we shall introduce a new proof of Small's formula
(Theorem~\ref{thm:small}), 
which is an analogue of the classical representation formula
for null curves in $\C^3$ (see the appendix).
We fix a Riemann surface $M^2$, which is not necessarily simply
connected.

Let
\begin{equation}\label{eq:null-curve}
   F = \begin{pmatrix} A & B \\ C & D \end{pmatrix}
\end{equation}
be a null curve in $\PSL(2,\C)$ defined on $M^2$.
\begin{definition}\label{def:null-gauss}
 For a non-constant null curve $F$ as in \eqref{eq:null-curve},
 we define
\[
   G : = \begin{cases}
	     \dfrac{dA}{dC}
               & \bigl(\text{if $(dA,dC)\not\equiv(0,0)$}\bigr),\\[8pt]
	     \dfrac{dB}{dD}
               & \bigl(\text{if $(dB,dD)\not\equiv(0,0)$}\bigr),
	  \end{cases}\quad
   g : = \begin{cases}
	     -\dfrac{dB}{dA}
               & \bigl(\text{if $(dA,dB)\not\equiv(0,0)$}\bigr),\\[8pt]
	     -\dfrac{dD}{dC}
               & \bigl(\text{if $(dC,dD)\not\equiv(0,0)$}\bigr).
	  \end{cases}
\]
 Since $F$ is null, 
\begin{align*}
   \hphantom{-}\frac{dA}{dC}&=\hphantom{-}\frac{dB}{dD}\qquad
   &\text{if $(dA,dC)\not\equiv (0,0)$ and  $(dB,dD)\not\equiv (0,0)$}\\
   -\frac{dB}{dA}&=-\frac{dD}{dC}\qquad
   &\text{if $(dA,dB)\not\equiv (0,0)$ and  $(dC,dD)\not\equiv (0,0)$}
\end{align*}
 hold.
 We call $G$ and $g$ the {\em hyperbolic Gauss map\/} and
 the {\em secondary Gauss map\/} of $F$, respectively.
\end{definition}
\begin{lemma}\label{lem:G-const}
 Let $F$ be a meromorphic null curve as in \eqref{eq:null-curve}.
 If either $dA\equiv dC\equiv 0$ or $dB\equiv dD\equiv 0$ holds,
 then the hyperbolic Gauss map $G$ is constant.
 Similarly, 
 if either $dA\equiv dB\equiv 0$ or $dC\equiv dD\equiv 0$ holds,
 then the secondary Gauss map $g$ is constant.
\end{lemma}
\begin{proof}
 Assume $dA\equiv dB\equiv 0$.
 Since $AD-BC=1$, we have
 \[
     0 = d(AD-BC) = DdA + A dD - BdC-C dB
                  = AdD - C dB.
 \]
 Here, since $(A,C)\not\equiv (0,0)$, 
 $dB/dD\in \C\cup\{\infty\}$ is constant.
 The other statements are proved in the same way.
\end{proof}
\begin{lemma}[\cite{UY2}, \cite{Yu}]\label{lem:g-bryant}
 Let $F$ be a non-constant null meromorphic curve such that
 the secondary Gauss map $g$ is non-constant.
 Set
\begin{equation}\label{eq:diffeq}
    F^{-1}dF = \alpha,\qquad
         \alpha =\begin{pmatrix}
		\alpha_{11} & \alpha_{12} \\
		\alpha_{21} & \alpha_{22}
	       \end{pmatrix}.
\end{equation}
 Then the secondary Gauss map $g$ of $F$ is represented as 
\[
    g = \frac{\alpha_{11}}{\alpha_{21}} = \frac{\alpha_{12}}{\alpha_{22}}.
\]
\end{lemma}
\begin{proof}
 Let $F$ be as in \eqref{eq:null-curve}.
 If $\alpha_{11}$ and $\alpha_{21}$ vanish identically,
 so is $\alpha_{22}$, because $\alpha$ is an $\Sl(2,\C)$-valued
 $1$-form.
 Then, since $dF=F\alpha$,
 we have $dA\equiv dB\equiv dD\equiv 0$, which implies
 $g$ is constant.
 Hence $(\alpha_{11},\alpha_{21})\not\equiv (0,0)$.
 Similarly, $(\alpha_{12},\alpha_{22})\not\equiv (0,0)$.
 Here $\det\alpha=0$ because $F$ is null.
 Hence we have $\alpha_{11}/\alpha_{21}=\alpha_{12}/\alpha_{22}$.
 
 Since $AD-BC=1$, it holds that $DdA-BdC=-AdD+CdB$.
 Then, using the relations $dB=-g\,dA$ and $dD=-g\,dC$, we have
 \[
   \frac{\alpha_{11}}{\alpha_{21}}
   = \frac{\hphantom{-}D dA-BdC}{-CdA+AdC}=\frac{-AdD+CdB}{-CdA+AdC}
   = \frac{g(-CdA+AdC)}{\hphantom{g(}-CdA+AdC\hphantom{)}}=g.
 \]
 This completes the proof.
\end{proof}
\begin{theorem}[Small \cite{S}]\label{thm:small}
 For an arbitrary pair of non-constant meromorphic functions
 $(G,g)$ on $M^2$ such that $g\not\equiv a\star G$ for any
 $a\in\PSL(2,\C)$,
 a meromorphic map $F$ given by \eqref{eq:small-formula} is a
 non-constant null curve in $\PSL(2,\C)$ whose hyperbolic Gauss map and
 secondary Gauss map are $G$ and $g$ respectively.

 Conversely, any meromorphic null curve in $\PSL(2,\C)$ whose hyperbolic
 Gauss map $G$ and secondary Gauss map $g$ are both non-constant are
 represented in this way.
\end{theorem}
An analogue of this formula for null curves in $\C^3$ is mentioned
in Appendix~\ref{app:A}.
\begin{proof}[Proof of Theorem~\ref{thm:small}]
 Let $(G,g)$ be a pair as in the statement of the theorem and 
 set as in \eqref{eq:small-formula}.
 Then
 \[
    \det F = -a\frac{db}{dG} + b \frac{da}{dG}
           = -a^2\frac{d}{dG}\left(\frac{b}{a}\right)
           = a^2\frac{dg}{dG} = 1,
 \]
 and 
 \[
    \frac{dF}{dG} =
        \begin{pmatrix}
	  \dfrac{dA}{dG} &
	  \dfrac{dB}{dG}\\[8pt]
	  \dfrac{dC}{dG} &
	  \dfrac{dD}{dG} 
	 \end{pmatrix}
       =\begin{pmatrix}
	   G\dfrac{d^2a}{dG^2} & G \dfrac{d^2 b}{dG^2} \\[8pt]
	   \hphantom{G}\dfrac{d^2a}{dG^2} & \hphantom{G} \dfrac{d^2 b}{dG^2} 
	 \end{pmatrix}.
 \]
 Hence $\rank dF\leq 1$, and $F$ is a meromorphic null curve  in
 $\PSL(2,\C)$.
 The hyperbolic Gauss map of $F$ is obtained as
 \[
    \frac{dA}{dC} 
      = \frac{dA/dG}{dC/dG} = G.
 \]
 On the other hand, the secondary Gauss map is obtained by 
 Lemma~\ref{lem:g-bryant} as
 \[
    \frac{\alpha_{11}}{\alpha_{21}}
       =\frac{\hphantom{-}D\dfrac{dA}{dG}-B\dfrac{dC}{dG}}{%
            -C\dfrac{dA}{dG}+A\dfrac{dC}{dG}}
       = -\frac{GD-B}{GC-A}
       = -\frac{b}{a}=g.
 \]

 Next, we prove that $F$ is non-constant.
 Assume $F$ is constant.
 Then by \eqref{eq:small-formula}, $da/dG=p=\mbox{constant}$.
 Thus we have $a = \sqrt{dG/dg}=p G + q$, where $p$ and $q$ are
 complex numbers.
 Hence
 \[
      \frac{dg}{dG} = \frac{1}{(pG+q)^2}.
 \]
 Integrating this, we have that $g$ is obtained as a
 M\"obius transformation of $G$, a contradiction.
 Thus  the first part of the theorem is proved.

 Conversely, let $F$ be a null curve as in \eqref{eq:null-curve}.
 By Definition~\ref{def:null-gauss}, we have
 \begin{equation}\label{eq:dAB}
  dA=G \, dC,\qquad dB=G \,dD.
 \end{equation}
 We set
 \[
   a:=GC-A,\qquad b:=GD-B.
 \]
 By \eqref{eq:dAB}, we have $da=C \,dG$ and $db=D \, dG$.
 Since $G$ is not constant, we have
 \begin{equation}\label{eq:CD}
  C=\frac{da}{dG},\qquad D=\frac{db}{dG}.
 \end{equation}
 Then $F$ can be expressed in terms of $a$ and $b$ as follows:
 \[
   F=\begin{pmatrix}
      G\displaystyle\frac{da}{dG}-a & 
      G\displaystyle\frac{db}{dG}-b \\[6pt]
      \hphantom{G}\displaystyle\frac{da}{dG}\hphantom{-a} & 
      \hphantom{G}\displaystyle\frac{db}{dG}\hphantom{-b}
     \end{pmatrix}.
 \]
 Since $\det F=1$, we have
 \begin{equation}\label{eq:det}
   \det
    \begin{pmatrix}
     -a & -b \\
     \displaystyle\frac{da}{dG} & \displaystyle\frac{db}{dG}
    \end{pmatrix}=1.
 \end{equation}
 Taking the derivative of this equation, 
 \begin{equation}\label{eq:det-D}
   \det
    \begin{pmatrix}
     -a & -b \\
     d\left(\displaystyle\frac{da}{dG}\right) & 
     d\left(\displaystyle\frac{db}{dG}\right)
    \end{pmatrix}=0
 \end{equation}
 holds.
 Here, since $g$ is non-constant, 
 $(dC,dD)\not\equiv 0$ by Lemma~\ref{lem:G-const}.
 Then by \eqref{eq:CD} and \eqref{eq:det-D}, it holds that
 \[
   g=\frac{dD}{dC}=-
     \frac{d(db/dG)}
          {d(da/dG)}=-\frac{b}{a}.
 \]
 This yields 
 \begin{equation}\label{eq:b-exp}
   b=-ga,\qquad 
   db=-(da) g-a(dg).
 \end{equation}
 Again by \eqref{eq:det}
 \[
   dG=\det
      \begin{pmatrix}
       -a & -b \\
       {da} & {db}
      \end{pmatrix}=
      \det
       \begin{pmatrix}
	-a & ag \\
	{da} & -(da) g-a(dg)
       \end{pmatrix}=a^2 dg.
 \]
 By this and \eqref{eq:b-exp}, we have
 $a=\sqrt{dG/dg}$ and $b=-ga$
 which implies \eqref{eq:small-formula}.
\end{proof}
By Theorem~\ref{thm:small}, we can prove the uniqueness of 
null curves with given hyperbolic Gauss map and secondary Gauss map.
Hence we have
\begin{corollary}\label{cor:null-moduli}
 Let $\N(M^2)$ be the set of non-constant null curves in $\PSL(2,\C)$ 
defined  on a Riemann surface 
$M^2$ with non-constant hyperbolic Gauss map and secondary Gauss map.
 Then $\N(M^2)$ corresponds bijectively to the set
 \[
  \left\{(G,g)\left|
     \begin{array}{l}
         \text{$G$ and $g$ are non-constant meromorphic functions on $M^2$}\\
         \text{such that $G\not\equiv a\star g$ for any $a\in\SL(2,\C)$.}
     \end{array}
  \right.\right\}.
 \]
\end{corollary}

It should be remarked that $(G,g)$ satisfies the following
important relation (see \cite{UY2}):
\begin{equation}
 S(g)-S(G)=2Q,
\end{equation}
where $Q$ is the Hopf differential of $F$ defined by
$Q:=(A dC-C dA)dg$ and $S$ is the Schwarzian derivative defined by
\[
  S(G)=
  \left[
  \left(\frac{G''}{G'}\right)'-\frac 12 \left(\frac{G''}{G'}\right)^2
  \right]\,dz^2 \qquad
  \left( '=\frac{d}{dz}\right)
\]
with respect to a local complex coordinate $z$ on $M^2$.
Though meromorphic $2$-differentials $S(g)$ and $S(G)$ depend on complex
coordinates, 
the difference $S(g)-S(G)$ does not depend on the choice of complex
coordinates.
\section{Legendrian curves in $\PSL(2,\C)$}
\label{sec:legendre}
In this section, we shall give a representation formula
for Legendrian curves in terms of 
two meromorphic functions $G$ and $G_*$,
which are called {\em the hyperbolic Gauss maps}.
We fix a Riemann surface $M^2$, which might not be simply
connected.
Let $E$ be a meromorphic Legendrian curve on $M^2$ as in
\eqref{eq:legendrian-curve}.
Since $AD-BC=1$, we can define two meromorphic functions
$G$ and $G_*$ as in \eqref{eq:legendrianG}.
We call $G$ and $G_*$ the {\em hyperbolic Gauss maps\/} of $E$.
(The geometric meaning of these hyperbolic Gauss maps is described in 
\cite{GMM1}.)
\begin{definition}\label{def:invariants-legendrian}
 Let $E$  be a  meromorphic Legendrian curve $E$ as in
 \eqref{eq:legendrian-curve}. 
 Then we can write 
 \begin{equation} \label{eq:omega-newdef}
  E^{-1}dE
   =\begin{pmatrix}
     0 & \theta \\
     \omega & 0
    \end{pmatrix},
 \end{equation}
 where $\omega$ and $\theta$ are meromorphic $1$-forms on $M^2$.
 We call $\omega$ the {\em canonical form\/}  and
 $\theta$ the {\em dual canonical form\/} of $E$.

 For a Legendrian curve $E$, we define another Legendrian curve $\hat E$
 by 
 \[
    \hat E = E 
        \begin{pmatrix} 
	 0 &  i \\
	 i & 0
	\end{pmatrix}.
 \]
 We call $\hat E$ the {\em dual\/} of $E$.
 The hyperbolic Gauss maps $\hat G$ and $\hat G_*$ of $\hat E$ satisfy 
  $\hat G = G_*$ and $\hat G_*= G$,
 and  the canonical form and the dual canonical form of $\hat E$
 are $\theta$ and $\omega$ respectively.
 Roughly speaking, the duality exchanges the role of $(G,\omega)$ and
 $(G_*,\theta)$.
\end{definition}
The following lemma holds.
\begin{lemma}\label{lem:can-form}
 For a non-constant meromorphic Legendrian curve $E$ as in
 \eqref{eq:legendrian-curve}, the following identities hold{\rm :}
 \begin{align}\label{eq:omega-def}
    \omega&=\begin{cases}
	     \dfrac{dA}{B}
             & 
	     \bigl(\text{if $dA\not\equiv0$ or $B\not\equiv 0$}\bigr),\\[8pt]
	     \dfrac{dC}{D}
	     &\bigl(\text{if $dC\not\equiv0$ or $D\not\equiv 0$}\bigr),\\[8pt]
	  \end{cases}\\
  \label{eq:theta-def}
    \theta&=\begin{cases}
	     \dfrac{dB}{A}
             & 
	     \bigl(\text{if $dB\not\equiv0$ or $A\not\equiv 0$}\bigr),\\[8pt]
	     \dfrac{dD}{C}
	     &\bigl(\text{if $dD\not\equiv0$ or $C\not\equiv 0$}\bigr).\\[8pt]
	  \end{cases}
 \end{align}
 Here $dA\not\equiv 0$ $($resp.~$B\not\equiv 0)$ means a $1$-form $dA$
 $($resp.~a function $B)$ is not identically $0$.
 In particular,  if all cases in \eqref{eq:omega-def} and \eqref{eq:theta-def}
 are well-defined,
 \[
     \omega= \frac{dA}{B}=\frac{dC}{D} \qquad\text{and}\qquad
     \theta= \frac{dB}{A}=\frac{dD}{C}
 \]
 hold.
\end{lemma}
\begin{proof}
 Since $E$ is Legendrian, $DdA-BdC=0$ holds, and 
 $\omega = A dC - C dA$ by \eqref{eq:omega-newdef}. 
 Hence we have
 \begin{align*}
    B \omega &= AB dC - BC dA = ADdA-BC dA = (AD-BC)dA = dA\\
\intertext{and}
    D \omega &= AD dC - CD dA = (AD-BC)dC = dC,
 \end{align*}
 which imply \eqref{eq:omega-def}.

 On the other hand, differentiating $AD-BC=1$, we have
 \[
     0 = d(AD-BC)=(DdA-BdC)+(AdD-CdB) = A dD - CdB.
 \]
 Since $\theta=DdB-BdD$, we have then
 \[
    A \theta = AD dB - AB dD = (AD-BC)dB = dB\qquad\text{and}\qquad
    C \theta = dD,
 \]
 which imply \eqref{eq:theta-def}.
\end{proof}
\begin{theorem}\label{thm:repG}
 Let $G$ and $G_*$ be non-constant meromorphic functions on $M^2$ such
 that $G$ is not identically equal to $G_*$.
 Assume  that
\begin{enumerate}
\renewcommand{\theenumi}{(\roman{enumi})}
\let\labelenumi\theenumi
 \item\label{item:xi-order}
      all poles of the $1$-form $\dfrac{dG}{G-G_*}$ are of order $1$,
      and
 \item\label{item:xi-form}
      $\displaystyle\int_{\gamma}\dfrac{dG}{G-G_*}\in \pi i \Z$ holds
      for each loop $\gamma$ on $M^2$.
\end{enumerate}
 Set
\begin{equation}\label{eq:xi}
    \xi(z) := c \exp\int_{z_0}^{z}\frac{dG}{G-G_*},
\end{equation}
 where $z_0\in M^2$ is a base point and $c\in \C\setminus\{0\}$ is
 an arbitrary constant.
 Then 
\begin{equation}\label{eq:repG2}
 E := \begin{pmatrix}
      G/\xi &      \xi G_{*}/(G-G_{*}) \\
      1/\xi &      \hphantom{G_*}\xi/(G-G_{*}) 
     \end{pmatrix}
\end{equation}
 is a non-constant meromorphic Legendrian curve in $\PSL(2,\C)$ whose
 hyperbolic Gauss maps are $G$ and $G_*$.
 The canonical form $\omega$ of $E$ is written as
\begin{equation}\label{eq:canG}
	 \omega=-dG/{\xi^2}.
\end{equation}
 Moreover, a point $p\in M^2$ is a pole of $E$ if and only if
 $G(p)=G_*(p)$ holds.

 Conversely, any meromorphic Legendrian curve in $\PSL(2,\C)$ with
 non-constant hyperbolic Gauss maps $G$ and $G_*$ is obtained 
 in this way.
\end{theorem}
\begin{proof}
 By the assumptions \ref{item:xi-order} and \ref{item:xi-form}, 
 $\xi^2$ is a meromorphic function on $M^2$.
 Hence $E$ as in \eqref{eq:repG2} is a meromorphic curve in $\PSL(2,\C)$.
 One can easily see that $\det E=1$ and  $DdA-BdC=0$, that is, $E$ is a
 Legendrian map with  hyperbolic Gauss maps $G$ and $G_*$.
 The canonical form $\omega$ is obtained as \eqref{eq:canG} using
 \[
     d\xi = \frac{\xi dG}{G-G_*}.
 \]
 Since $G=A/C$ are non-constant, so is $E$.

 Next, we fix a point $p\in M^2$.
 By a matrix multiplication $E\mapsto \widetilde E=aE$
 ($a\in\SL(2,\C)$), 
 we have another Legendrian map $\widetilde E$ 
 with hyperbolic Gauss maps $\widetilde G=a\star G$ and
 $\widetilde G_*=a\star G_*$, where $\star$ denotes the M\"obius
 transformation \eqref{eq:moebius}.
 If necessary replacing $E$ by $\tilde E$, we may assume $G(p)\neq\infty$ and 
 $G_*(p)\neq\infty$.
 Let $z$ be a local complex coordinate on $M^2$ such that  $z(p)=0$.

 Assume $E$ is holomorphic at $p$.
 Then by \eqref{eq:repG2}, 
 $CD = 1/(G-G_*)$ is holomorphic at $p$.
 Hence we have $G(p)\neq G_*(p)$.
 On the other hand, if $G(p)\neq G_*(p)$, $\xi$ is holomorphic
 at $p$ and $\xi(p)\neq 0$.
 Then by \eqref{eq:repG2}, $E$ is holomorphic at $p$.
 Thus, we have shown that $\{p\in M^2|G(p)\neq G_*(p)\}$
 is the set of poles of $E$.

 Finally, we shall prove the converse statement.
 Let $E$ as in \eqref{eq:legendrian-curve} be a meromorphic Legendrian
 curve. 
 Then by \eqref{eq:omega-def}, we have
 \begin{equation} \label{eq:dG}
    dG=d\left(\frac AC \right)=\frac{CdA-AdC}{C^2}=
           -\frac{\omega}{C^2}=-\frac{dC}{C^2D}.
 \end{equation}
 On the other hand, we have
 \begin{equation}\label{eq:G_G*}
  G-G_*=\frac{A}{C}-\frac{B}{D}
   =\frac{AD-BC}{CD}=\frac1{CD}.
 \end{equation}
 By \eqref{eq:dG} and \eqref{eq:G_G*} 
 \begin{equation}\label{eq:log-c}
     d\log C = -\frac{dG}{G-G_{*}}
 \end{equation}
 holds.
 Since $E$ is a meromorphic map into $\PSL(2,\C)$, 
 $C$ is written as in the form
 $\sqrt{h}\hat C$, where $h$ and $\hat C$ are meromorphic 
 functions.
 Then if we set $\xi$ as in \eqref{eq:xi}, $\xi^2$ is a meromorphic 
 function on $M^2$.
 Hence we have \ref{item:xi-order} and \ref{item:xi-form} 
 in the statement of the theorem.
 Integrating \eqref{eq:log-c}, we have $C=1/\xi$ and $A=GC=G/\xi$. 
 Moreover, since
 \[
     1 = AD-BC = \left(\frac{G}{\xi}\right)D - G_*D\left(\frac{1}{\xi}\right)
       = \frac{D}{\xi} (G-G_*),
 \]
 we have $D=\xi/(G-G_*)$ and $B=G_*D=G_*\xi/(G-G_*)$.
 Thus we obtain \eqref{eq:repG2}.
\end{proof}

As a corollary of Theorem~\ref{thm:repG}, 
we give a Small-type formula for Legendrian curves,
which has appeared implicitly in \cite{GMM1} by a different
method.
\begin{corollary}\label{cor:legendre}
 For an arbitrary pair $(G,\omega)$ of a non-constant
 meromorphic function and a non-zero meromorphic $1$-form on $M^2$,
 a meromorphic map 
 \begin{equation}\label{eq:sm-analogue2}
  E=\begin{pmatrix}
     A & dA/\omega \\
     C & dC/\omega 
    \end{pmatrix}\qquad
    \left(C:=i\sqrt{\frac{\omega}{dG}},~
          A:=GC\right)
 \end{equation}
 is a meromorphic Legendrian curve in $\PSL(2,\C)$ whose hyperbolic
 Gauss map and canonical form are $G$ and $\omega$, respectively.

 Conversely, let $E$ be a meromorphic Legendrian curve in $\PSL(2,\C)$
 defined on $M^2$ with the non-constant hyperbolic Gauss map $G$ and
 the non-zero canonical form $\omega$.
 Then $E$ is written as in \eqref{eq:sm-analogue2}.
\end{corollary}
\begin{remark*}
  There is a correponding simple formula (with no integration)
  for Legendrian curves  in $\C^3$ as follows:
  A meromorphic map $E\colon{}M^2\to \C^3$ is called {\em Legendrian}
  if the pull-back of the holomorphic
  contact form $dx^1-x^3\,dx^2$ vanishes,
  where $(x^1,x^2,x^3)$ is the canonical coordinate system on $\C^3$.
  For a pair $(f,g)$ of meromorphic functions on a Riemann surface
  $M^2$, $E:=(f,g,df/dg)$ trivially gives a meromorphic Legendrian
  curve, which is an analogue of \eqref{eq:sm-analogue2}.
\end{remark*}
\begin{proof}[Proof of Corollary~\ref{cor:legendre}]
 If we set $E$ by \eqref{eq:sm-analogue2},
 we have $AD-BC=1$ and $DdA-BdC=0$.
 Hence $E$ is a meromorphic Legendrian map.
 
 Conversely,
 let $E$ be a meromorphic Legendrian  curve on $M^2$ 
 with the non-constant hyperbolic Gauss map $G$ and
 the non-zero canonical form $\omega$.
 Then by \eqref{eq:dG}, we have
 \[
   C=i\sqrt{\frac{\omega}{dG}},\qquad
   A=GC.
 \]
 On the other hand, by Lemma~\ref{lem:can-form}, we have
 $B=dA/\omega$ and $D=dC/\omega$.
 Hence we have \eqref{eq:sm-analogue2}.
\end{proof}

We have the following corollary:
\begin{corollary}\label{cor:legendrian-moduli}
 Let $\L(M^2)$ be the set of meromorphic Legendrian curves in \newline
 $\PSL(2,\C)$  defined on a Riemann surface $M^2$ with 
 non-constant hyperbolic Gauss maps and non-zero canonical forms.
 Then $\L(M^2)$ corresponds 
 bijectively to the following set{\rm :}
 \[
    \left\{(G,\omega)\,\left|
      \begin{array}{l}
       \text{$G$ is a non-constant meromorphic function on $M^2$,}\\
       \text{and $\omega$ is a non-zero meromorphic $1$-form on $M^2$}.
      \end{array}
    \right.\right\}.
 \]
\end{corollary}

The symmetric product of the canonical form $\omega$ and the 
dual form $\theta$ 
\begin{equation}\label{eq:legendre-hopf}
 Q:=\omega\theta
\end{equation}
is called the {\em Hopf differential\/} of the Legendrian curve.
By \eqref{eq:dG}, we have
\[
   dG=-\frac{\omega}{C^2}.
\]
Similarly, it holds that
\[
   dG_*=\frac{\theta}{D^2}.
\]
Thus, by \eqref{eq:G_G*} we have
\[
     Q=-C^2 D^2 dG\, dG_*=-\frac{dG\, dG_*}{(G-G_*)^2}.
\]
As pointed out in \cite{GMM1}, the following identities
hold:
\[
   S(g)-S(G)=2Q,\qquad S(g_*)-S(G_*)=2Q,
\]
where $g$ (resp. $g_*$) is a meromorphic function defined on
the universal cover of $M^2$ such that $dg=\omega$ 
(resp.~$dg_*=\theta$).

\section{Examples of flat surfaces in $H^3$}
\label{sec:example}
As an application of Corollary~\ref{cor:legendre}, we 
shall give new examples of flat surfaces in hyperbolic $3$-space $H^3$.
Though these examples might have singularities,
All of them are obtained as projections of Legendrian immersions
into the unit cotangent bundle $T_1^*H^3$.
Usually, a projection of a Legendrian immersion is called a 
(wave) front.
So we call them {\em flat fronts}.
For details, see \cite{KUY}.

Hyperbolic 3-space $H^3$ has an expression 
\[
  H^3=\PSL(2,\C)/\PSU(2)=\{aa^{*}\,;\,a\in\PSL(2,\C)\}\qquad
  (a^*={}^t\bar a).
\]
As shown in \cite{GMM1}, the projection 
\[
  f:=EE^*\colon{}M^2\longrightarrow H^3
\]
of a holomorphic Legendrian curve $E\colon{}M^2\to \PSL(2,\C)$ is a flat
immersion if $f$ induces positive definite metric on $M^2$.
For a Legendrian curve $E$, we can write
\begin{equation}\label{eq:gmm-formula}
    E^{-1}dE = \begin{pmatrix}
		 0 & \theta \\
		 \omega & 0 
		\end{pmatrix}.
\end{equation}
Then the first fundamental form $ds^2$
and second fundamental form $d\sigma^2$ of $f$ is written as
\begin{align}\label{eq:fund-form}
   ds^2&= \omega\theta + \overline{\omega\theta}
            +|\omega|^2 + |\theta|^2
       = (\omega+\bar\theta)(\bar\omega+\theta),\\
 \label{eq:fund-form-2}
   d\sigma^2&=\pm(|\omega|^2-|\theta|^2).
\end{align}

\begin{figure}
\small
\begin{center}
\begin{tabular}{c@{\hspace{1cm}}c@{\hspace{1cm}}c}
 \includegraphics[width=3cm]{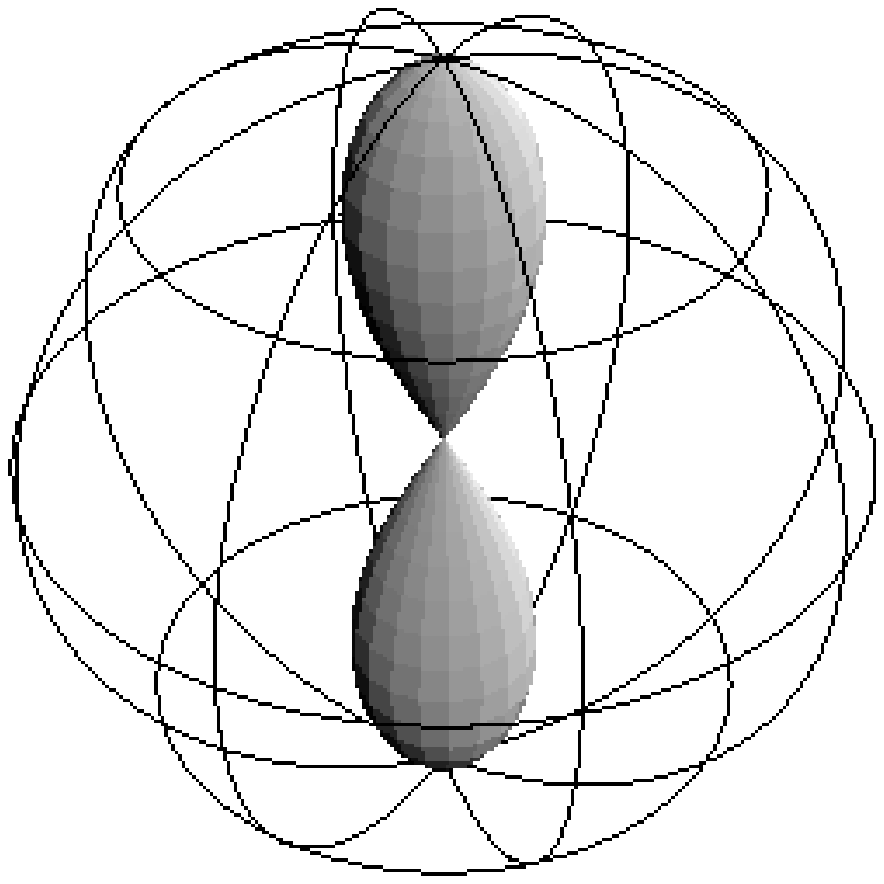} &
 \includegraphics[width=3cm]{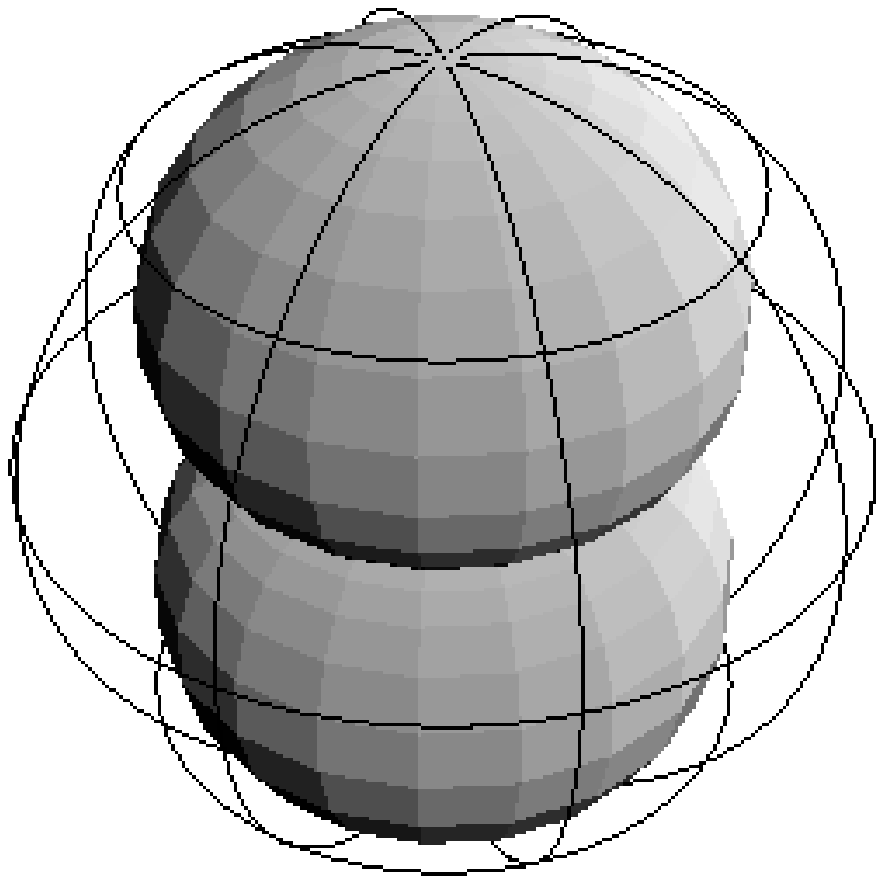} &
 \includegraphics[width=3cm]{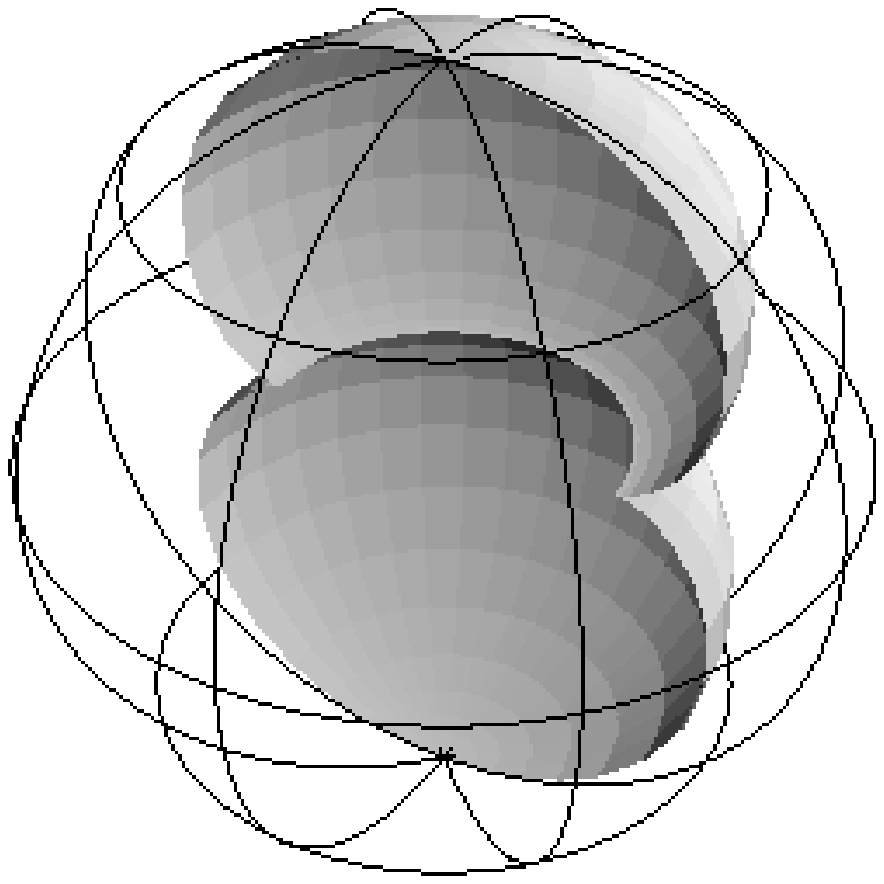} \\
 $k=0.5$ &
 $k=1.5$ &
 $k=1.5$ (half cut)
\end{tabular}
\end{center}
\caption{Flat fronts of revolution (Example~\ref{exa:revolution}).}
\label{fig:revolution}
\end{figure}

\begin{figure}
 \footnotesize
\begin{center}
  \begin{tabular}{c@{\hspace{2cm}}c}
    \includegraphics[width=3cm]{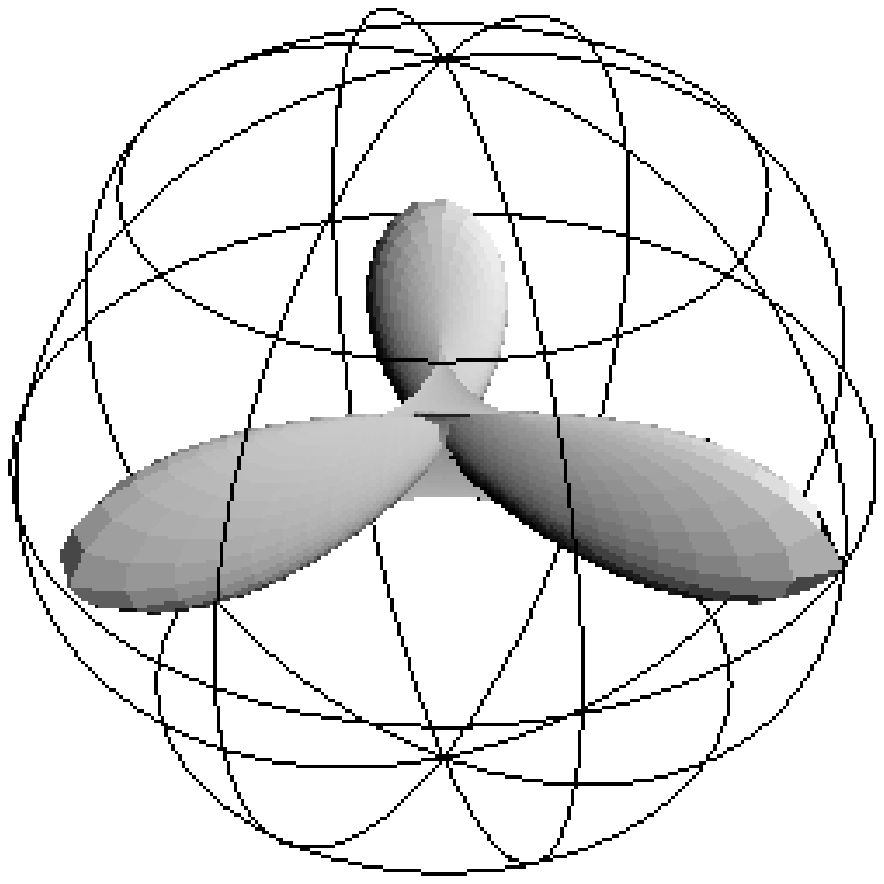} &
    \includegraphics[width=3cm]{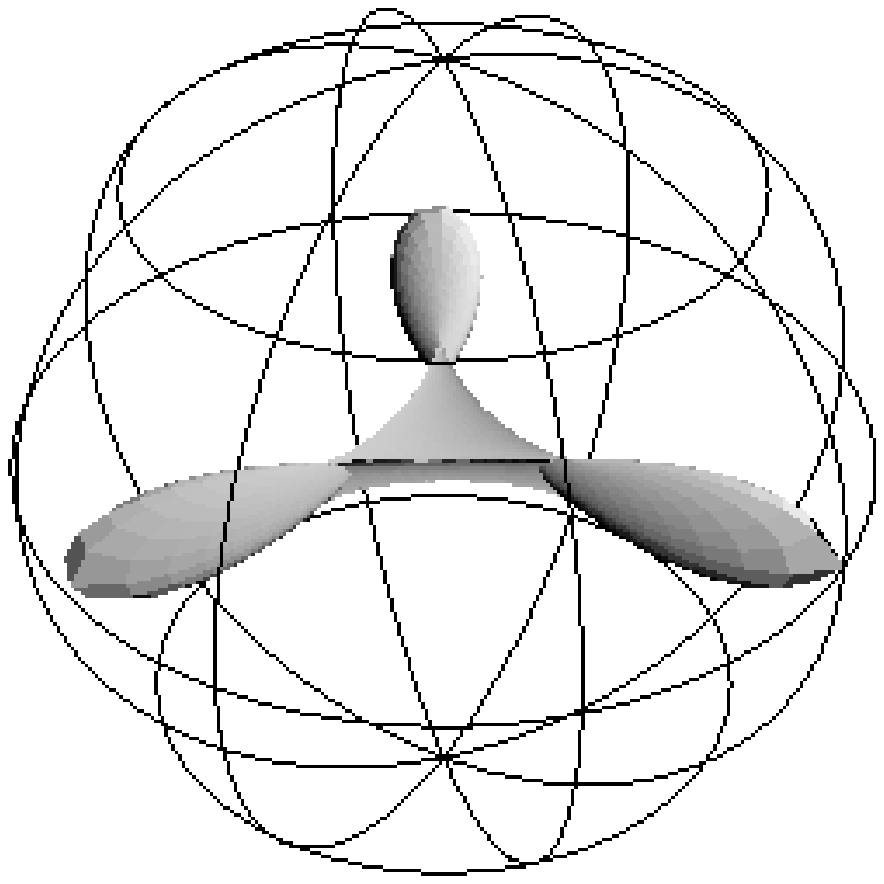}  \\
    $k=0.8$ &
    $k=1.0$ \\[4ex]
    \includegraphics[width=3cm]{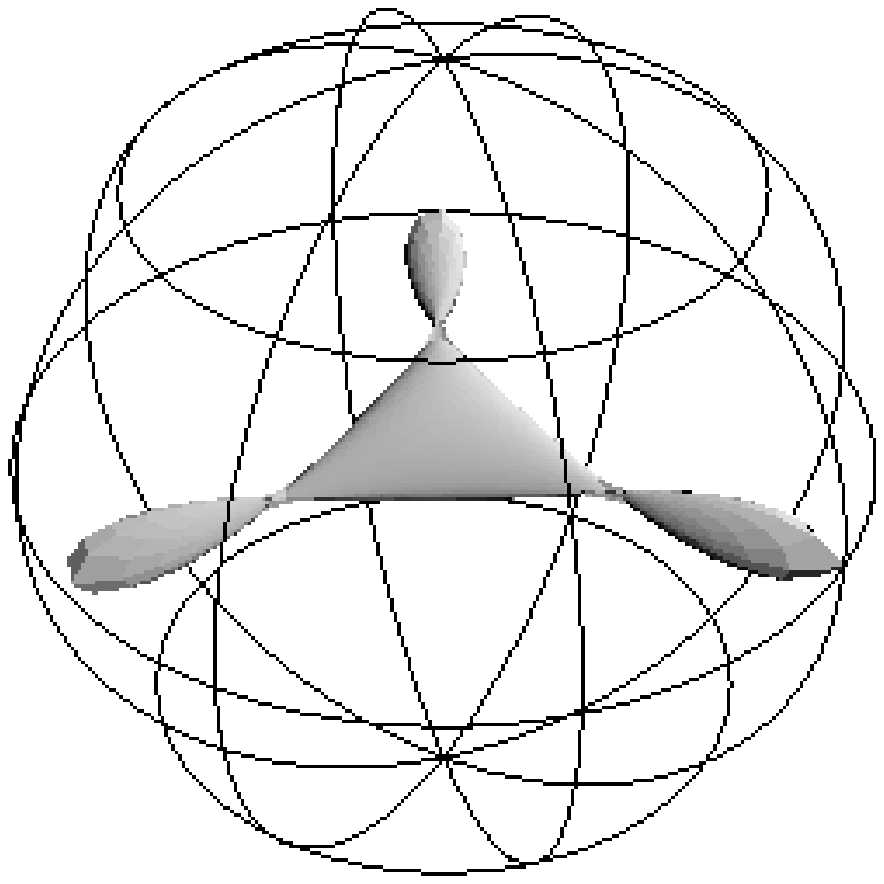} &
    \includegraphics[width=3cm]{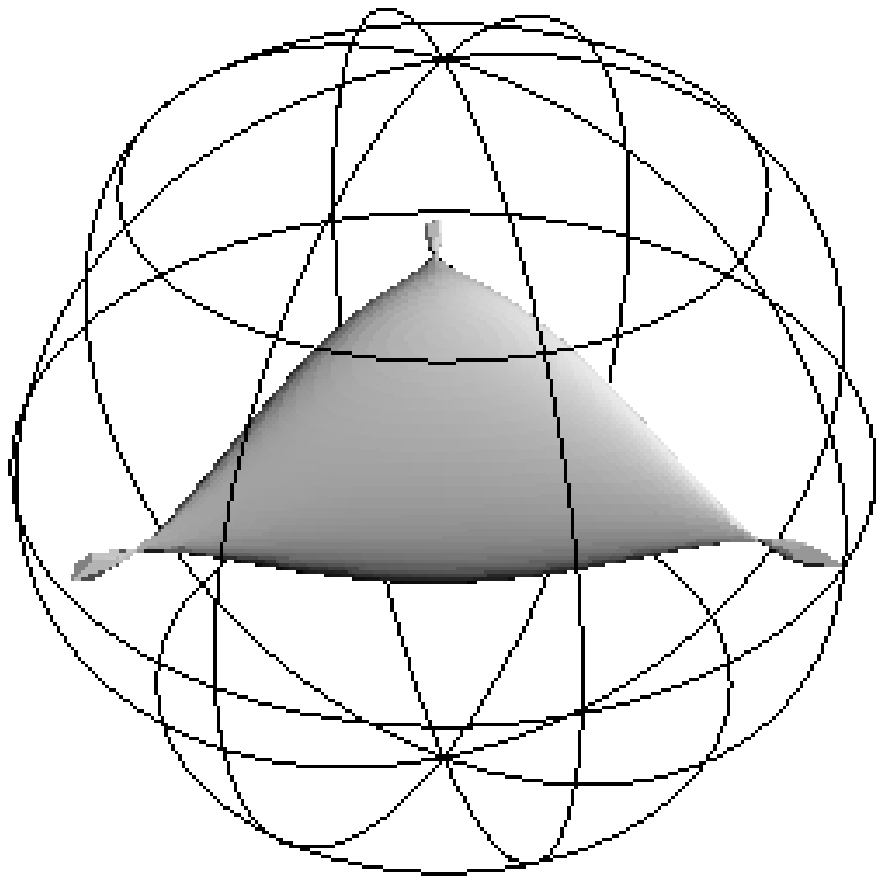}\\
    $k=1.2$ &
    $k=1.8$ 
  \end{tabular}
\end{center}
\caption{Parallel family of flat fronts in
 Example~\ref{ex:jorge-meeks} for $n=3$}
\label{fig:jorge-meeks}
\end{figure}

Common zeros of $\omega$ and $\theta$
correspond to  branch points of the surface where
the first fundamental form vanishes.
At the point where $|\omega|=|\theta|$, 
$ds^2$ in \eqref{eq:fund-form} is written as
\[
   ds^2=\frac{\omega}{\theta}(\bar\omega+\theta)^2
\]
which implies the metric degenerates at these points.
Let $\nu$ be the unit normal vector field of $f$.
For each $p\in M^2$, the asymptotic class of the geodesic with initial
velocity $\nu(p)$ (resp.~$-\nu(p)$) determines a point $G(p)$
(resp.~$G_{*}(p)$) of the ideal boundary of $H^3$ which is 
identified with $\C\cup\{\infty\}=\CP^1$.
Then $G$ and $G_*$ coincide with the hyperbolic Gauss maps of the lift $E$.
\begin{example}[Surfaces equidistant from a geodesic]
\label{exa:cylinder}
 Let $M^2=\C\setminus\{0\}$ and 
 \[
     G = z, \qquad \omega = \frac{k}{2z}\,dz\qquad
     (k>0).
 \]
 Then by Corollary~\ref{cor:legendre}, the corresponding Legendrian curve
 $E$ is written by
 \[
     E = \frac{i}{\sqrt{2}}
         \begin{pmatrix}
	   \sqrt{kz} & \sqrt{\dfrac{z}{k}}\\[6pt]
	   \sqrt{\dfrac{k}{z}} & -\dfrac{1}{\sqrt{kz}}
	  \end{pmatrix}.
 \]
 Then corresponding flat surface $f=EE^*$ is a surface equidistant
 from a geodesic in $H^3$ .
 The hyperbolic Gauss maps of $f $ are given by $(G,G_*)=(z,-z)$
 (see Figure~\ref{fig:cylinder}, and see also \cite[page 426]{GMM1}).
\begin{figure}
\begin{center}
\small
 \begin{tabular}{c@{\hspace{2cm}}c}
  \includegraphics[width=3cm]{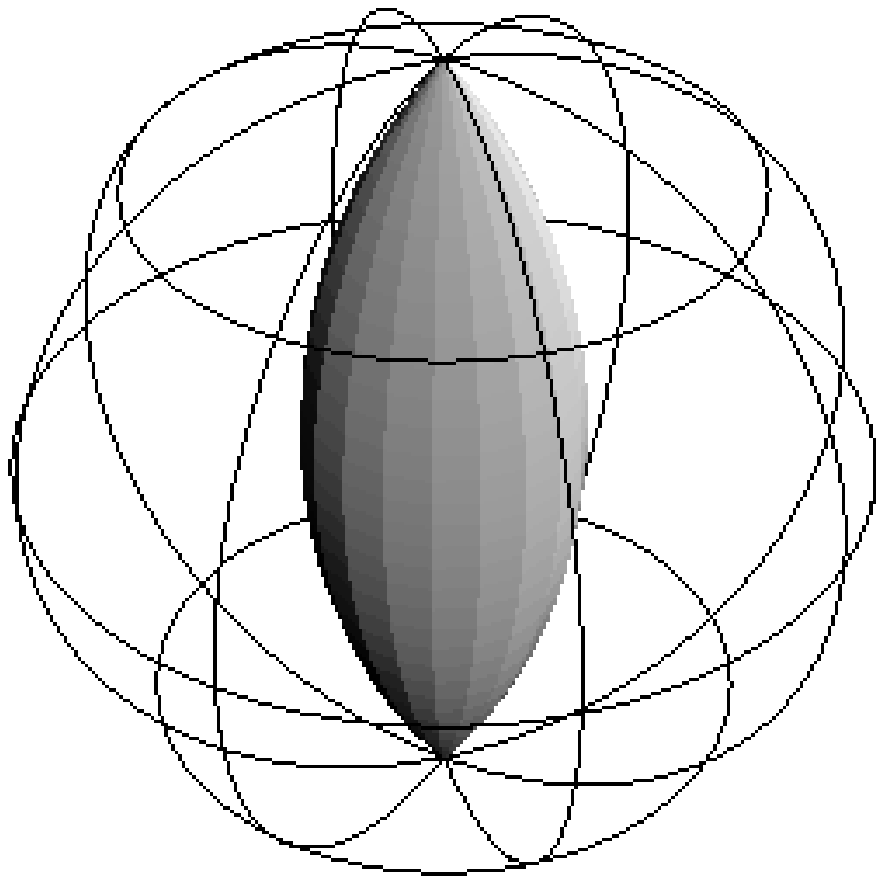} &
  \includegraphics[width=3cm]{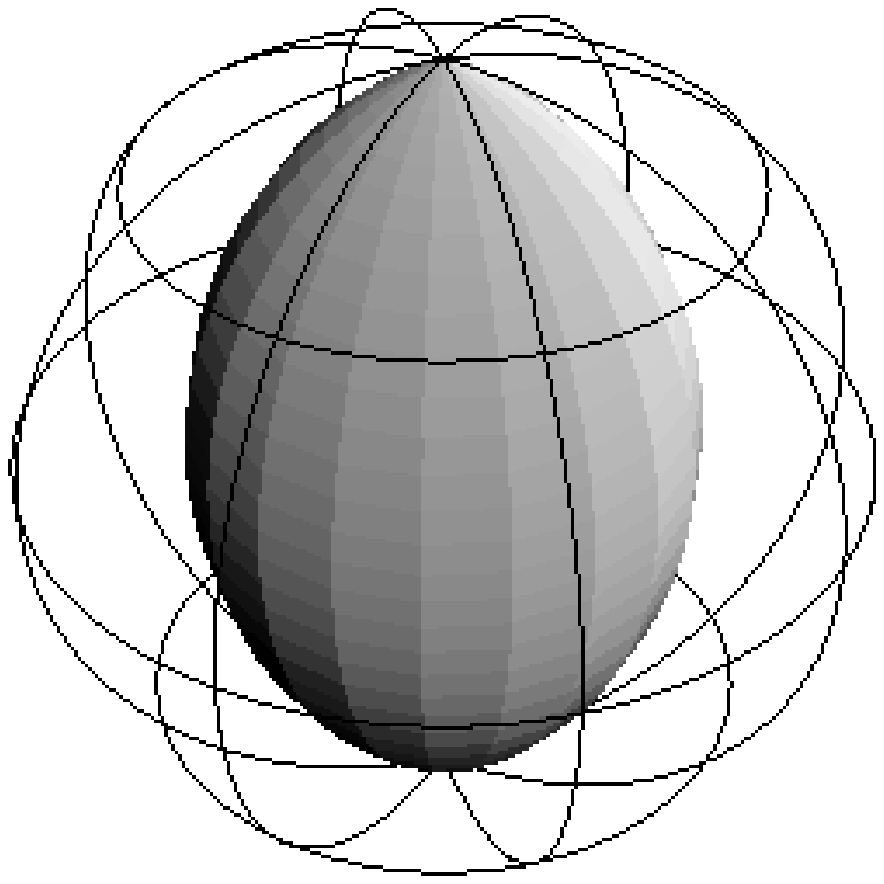} \\
  $k=\sqrt{2}$ & $k=2$
 \end{tabular}
\end{center}
\caption{Surfaces equidistant from a geodesic
 (Example~\ref{exa:cylinder}).
 The figures are shown in the Poincar\'e model of $H^3$.}
\label{fig:cylinder}
\end{figure}

\end{example}
\begin{example}[Flat fronts of revolution]
\label{exa:revolution}
 Let $M^2=\C\setminus\{0\}$
 and set
 \[
    G = \sqrt{\frac{\mu-1}{\mu+1}}z\qquad 
    \text{and} \qquad
    \omega = \frac{\sqrt{1-\mu^2}}{2}z^{\mu-1}\,dz 
    \qquad (\mu\in\R_+\setminus\{1\}).
 \]
 If $\mu\not\in\Z$, $\omega$ is not well-defined on $M^2$, but defined
 on the universal cover $\widetilde M^2$ of $M^2$.
 If we consider $G$ as a function on $\widetilde M^2$, 
 the corresponding Legendrian curve $E\colon{}\widetilde M^2\to\PSL(2,\C)$
 is given as
 \[
    E = \frac{i}{\sqrt{2}}
        \begin{pmatrix}
	 z^{(\mu+1)/2} & (\mu+1) z^{-(\mu-1)/2}\\
	 z^{(\mu-1)/2} & (\mu-1) z^{-(\mu+1)/2}
        \end{pmatrix}.
 \]
 Let $\tau$ be  the deck transformation of $\widetilde M^2$
 corresponding to the loop on $M^2$ surrounding $0$.
 Then
 \[
     E\circ\tau = E
                  \begin{pmatrix}
		   -e^{\pi i \mu} & 0 \\
		      0           & -e^{-\pi i \mu}
		  \end{pmatrix}
 \]
 holds.
 Hence the corresponding surface $f=EE^*$ is well-defined on $M^2$. 
 The dual canonical form $\theta$ as in \eqref{eq:gmm-formula} is given by
 \[
     \theta = \frac{\sqrt{1-\mu^2}}{2}z^{-\mu-1}\,dz.
 \]
 Then the metric induced by $f$ degenerates on the set $\{|z|=1\}$ when
 $\mu\neq 0$
 (see Figure~\ref{fig:revolution}).
 The hyperbolic Gauss maps of $f$ are given by
 \[
     (G,G_*)=\left(\sqrt{\frac{\mu-1}{\mu+1}}z,
                   i\sqrt{\frac{\mu+1}{\mu-1}}z\right).
 \]
\end{example}
\begin{example}[Flat fronts with dihedral symmetry]
\label{ex:jorge-meeks}
 Let $n\geq 2$ be an integer.
 We set 
\[
    M^2:=\C\cup\{\infty\}\setminus\{1,\zeta,\dots,\zeta^{n-1}\}
    \qquad\left(\zeta = \exp\frac{2\pi i}{n}\right).
\]
 and let $\pi\colon{}\widetilde M^2\to M^2$ be the universal cover of
 $M^2$.
 Let 
\begin{equation}\label{eq:jm-data}
   G_0(z)=z\quad\text{and}\quad 
   \omega = k(z^n-1)^{-2/n}\,dz\qquad (k>0),
\end{equation} 
 where $z$ is the canonical coordinate on $\C$.
 Then $G:=G_0\circ\pi$ and $\omega$ are considered as a meromorphic
 function and a holomorphic $1$-form on $\widetilde M^2$.
 Then by Corollary~\ref{cor:legendre}, there exists a holomorphic Legendrian 
 curve  $E\colon{}\widetilde M^2\to\PSL(2,\C)$.
 Let $\tau_j$ be a deck transformation of 
 $\pi\colon{}\widetilde M^2 \to M^2$
 corresponding a loop on $M^2$ around $\zeta^j$ ($j=0,\dots,n-1$).
 Then we have 
 \[
     G\circ\tau_j=G,\qquad
     \omega\circ\tau_j=\zeta^{-2}\omega.
 \]
 Hence by \eqref{eq:sm-analogue2}, we have
 \[
     E\circ\tau_j = E\begin{pmatrix}
		           \zeta^{-1} & 0 \\
		           0         & \zeta
		      \end{pmatrix} \qquad (j=0,\dots,n-1).
 \]
 This implies $f:=EE^*$ is well-defined on $M^2$ itself.
 Thus, we have a one parameter family of flat surfaces in $H^3$,
 parametrized by $k$ in \eqref{eq:jm-data}.
 The parameter $k$ corresponds to a parallel family of flat surfaces
 (see \cite[page 426]{GMM1}).
 Moreover, by \eqref{eq:fund-form}, one can see that each end $\zeta^j$
 is complete.
 On the other hand, at the points where $|\omega|=|\theta|$, 
 the immersion $f$ has singularities.
 The automorphisms of $M^2$ as
 \[
      z\longmapsto \zeta z\,, \qquad z\longmapsto 1/z
 \]
 do not change the first and second fundamental forms as in
 \eqref{eq:fund-form}.
 This implies such  surfaces have dihedral symmetry (see
 Figure~\ref{fig:jorge-meeks}).
 The hyperbolic Gauss maps of $f$ are given by
 \[
    (G,G_*)=\left(z, z^{1-n}\right).
 \]
\end{example}
\begin{example}[A flat front with tetrahedral symmetry]%
\label{ex:tetra}
 Let 
\[
    M^2 = \C\cup\{\infty\}\setminus\{1,\zeta,\zeta^2,\infty\}\qquad
          \left(\zeta=\exp\frac{2\pi i}{3}\right)
\]
 with $\pi\colon{}\widetilde M^2\to M^2$ the universal covering.
 Set 
\[
    G = G_0\circ\pi,\qquad G_0(z)=z
    \quad\text{and}\quad
    \omega =k(z^3-1)^{-1/2}\,dz \quad (k>0).
\]
 Then, in the same way as in Example~\ref{ex:jorge-meeks}, 
 we have a one parameter family of flat surfaces with four complete
 ends at $z=1,\zeta,\zeta^2,\infty$.
 Such surfaces have the tetrahedral symmetry.
 The hyperbolic Gauss maps of $f$ are given by
 \[
    (G,G_*) =\left( z ,  \frac{4-z^3}{3z^2}\right).
 \]
\end{example}
 In  Figures~\ref{fig:revolution} and \ref{fig:jorge-meeks}, 
 it seems that the surfaces admit singularities.
 It might be interesting problem to study singularities of flat
 fronts (see \cite{KUY}).

\appendix
\section{Null curves in $\C^3$}
\label{app:A}
Let $M^2$ be a Riemann surface, which is not necessarily simply
connected.
A meromorphic map  $F=(F^1,F^2,F^3)\colon M^2 \to \C^3$ 
is said to be {\it null\/} 
if the $\C^3$-valued $1$-form $dF$ is null, that is,  
\begin{equation}\label{2-1}
   \sum_{j=1}^3 dF^j\cdot dF^j =0. 
\end{equation}
It is well-known that a minimal surface in $\R^3$ 
is locally given by the projection of a null curve in $\C^3$ 
to $\R^3$.  

\medskip

For a null meromorphic map $F=(F^1,F^2,F^3)$, we put  
\begin{equation}\label{2-2}
  \omega:=d(F^1-i F^2),\qquad
         g:=dF^3/\omega. 
\end{equation}
Then we have 
\begin{equation}\label{2-3}
     dF=\frac 12\bigl((1-g^2)\omega,
                      i(1+g^2)\omega,
                       2g\omega\bigr)
\end{equation}
by \eqref{2-1}. Conversely, by integrating \eqref{2-3} 
for a given pair $(g,\omega)$, we obtain a null meromorphic map $F$. 
The integration of \eqref{2-3} is known as the 
{\it Weierstrass formula\/} and the pair 
$(g,\omega)$ is called the {\it Weierstrass data\/} of $F$. 

On the other hand, let $F \colon M^2 \to \C^3$ be 
a meromorphic map defined by 
\begin{equation}\label{eq:dar}
  F=\begin{pmatrix} 
             1 & g & \hphantom{i}(1-g^2)/2 \\
             i & ig& i(1+g^2)/2 \\
             0 & -1   &     g      
    \end{pmatrix}
    \begin{pmatrix} h\hphantom{_1} \\ h_1\\ h_2\end{pmatrix}
    \qquad \left(
              h_1=\frac{dh}{dg},~ h_2 = \frac{dh_1}{dg}
           \right)
\end{equation}
for a pair $(g,h)$ of two meromorphic functions, 
then $F$ is null. 
Conversely, any null meromorphic map $F \colon M^2 \to \C^3$ 
is represented by this formula \eqref{eq:dar}.
The Weierstrass formula \eqref{2-3} and the formula \eqref{eq:dar} are
related by $(g,\omega)=(g,dh_2)$. 

The remarkable feature of the formula \eqref{eq:dar} is 
that arbitrary null meromorphic maps 
can be represented in the integral-free form. 

\medskip

We introduce here a way to derive the formula \eqref{eq:dar}. 
\smallskip

Let $F \colon M^2 \to \C^3$ be a null curve and $(g, \omega)$ its 
Weierstrass data. We let  
\begin{equation}\label{null1}
  h_2:=F^1-iF^2,\quad
  \psi:=-F^1-iF^2,\quad
  \varphi:=F^3, 
\end{equation}
then their differentials satisfy
\begin{align}
\label{s1} 
  & dh_2=\omega \\
\label{s2}
  & d\varphi=g\omega\\
\label{s3}
  & d\psi=g^2 \omega. 
\end{align}
Now, we define a function $h_1$ by 
\begin{equation}\label{null2}
  \varphi=h_2g-h_1.
\end{equation}
Using \eqref{s1}  and \eqref{s2}, we compute that 
\[
   g\omega=d\varphi=d(h_2g-h_1)=g\omega+h_2dg-dh_1, 
\]
hence
\begin{equation} \label{null3}
   h_2=dh_1/dg. 
\end{equation}
Moreover, we define a function $h$ by 
\begin{equation}\label{null4}
    \psi=h_2g^2-2h_1g+2h,  
\end{equation}
then 
\begin{align*}
   g^2\omega&=d(h_2g^2-2h_1g+2h)\\
            &=g^2dh_2+2h_2g\,dg-2g\,dh_1-2h_1\,dg+2dh\\
            &=g^2\omega+2h_2g\,dg-2h_2g\,dg-2h_1dg+2dh\\
            &=g^2\omega -2h_1\,dg+2dh, 
\end{align*}
by \eqref{s1}--\eqref{s3}, hence 
\begin{equation} \label{null5}
   h_1=dh/dg. 
\end{equation}
Substituting \eqref{null2}--\eqref{null5} into \eqref{null1}, 
we obtain the formula \eqref{eq:dar}.

\end{document}